\documentclass[12pt,reqno]{amsart}
\usepackage{amssymb,amscd,amsbsy}
\newcommand{\lb}{\linebreak}

\newcommand{\z}{\zeta}

\newcommand{\vt}{\vartheta}

\newcommand{\U}{\Upsilon}

\newcommand{\B}{{\mathcal B}}

\newcommand{\h}{{\mathcal H}}

\newcommand{\C}{{\Bbb C}}
\newcommand{\T}{{\Bbb T}}
\newcommand{\pp}{{\Bbb P}}
\newcommand{\dd}{{\Bbb D}}

\newcommand{\mm}{{\Bbb M}}
\newcommand{\0}{{\boldsymbol{0}}}

\newcommand{\bs}{\boldsymbol}

\newcommand{\m}{{\boldsymbol m}}
\newcommand{\bS}{{\boldsymbol S}}

\newcommand{\rf}[1]{(\ref{#1})}

\newcommand{\df}{\stackrel{\mathrm{def}}{=}}
\newcommand{\dist}{\operatorname{dist}}

\newcommand{\trace}{\operatorname{trace}}
\newcommand{\rank}{\operatorname{rank}}

\newcommand{\eeq}{\end{equation}}
\newcommand{\beq}{\begin{equation}}
\newcommand{\bay}{\begin{eqnarray}}
\newcommand{\ba}{\begin{align*}}
\newcommand{\ea}{\end{align*}}
\newcommand{\ey}{\end{eqnarray}}
\newcommand{\bey}{\begin{eqnarray*}}
\newcommand{\eey}{\end{eqnarray*}}

\newcommand{\imp}{\Rightarrow}
\newcommand{\be}{\infty}

\newcommand{\bl}{\blacksquare}
\newcommand{\ess}{\operatorname{ess}}

\newcommand{\Pf}{{\bf Proof. }}

\newcommand{\ov}{\overline}

\newtheorem{thm}{\hspace{\parindent}Theorem}[section]

\begin{document}

\newcommand{\vse}{\vspace{.2in}}
\numberwithin{equation}{section}

\title{\bf Analytic approximation of matrix functions\\ and dual extremal functions}
\author{V.V. Peller}
\thanks{The author is partially supported by NSF grant DMS 0700995}

\maketitle

\begin{abstract}
We study the question of the existence of a dual extremal function for a bounded
matrix function on the unit circle in connection with the problem of approximation by analytic matrix functions. We characterize the class of matrix functions, for which a dual extremal function exists
in terms of the existence of a maximizing vector of the corresponding Hankel operator and
in terms of certain special factorizations that involve thematic matrix functions.
\end{abstract}

\section{\bf Introduction}
\setcounter{equation}{0}

\

In this paper we consider the problem of approximation of bounded
matrix-valued functions on the unit circle $\T$ by bounded analytic 
matrix functions in the unit disk $\dd$.  In other words, for $\Phi\in L^\be(\mm_{m,n})$
(i.e., $\Phi$ is a bounded function that takes values in the space $\mm_{m,n}$ of $m\times n$ matrices),
we search for a matrix function $F\in H^\be(\mm_{m,n})$ (i.e., $F$ is a bounded analytic function in 
$\dd$ with values in $\mm_{m,n}$) such that
\bay
\label{dis}
\|\Phi-F\|_{L^\be}=\dist_{L^\be}\big(\Phi,H^\be(\mm_{m,n})\big).
\ey
Here for a function $G\in L^\be(\mm_{m,n})$,
$$
\|G\|_{L^\be}\df\ess\sup_{\z\in\T}\|G(\z)\|_{\mm_{m,n}},
$$
where for a matrix $A\in\mm_{m,n}$, the norm $\|A\|_{\mm_{m,n}}$ is the norm of $A$
as an operator from $\C^n$ to $\C^m$. It is well known (and follows easily from a compactness argument) that the distance on the right-hand side of \rf{dis} is attained. A matrix function $\Phi$ is called {\it badly approximable} if the zero matrix function is a best approximant to $\Phi$ or, in other words,
$$
\|\Phi\|_{L^\be}\ge\|\Phi-F\|_{L^\be}\quad\mbox{for every}\quad F\in H^\be(\mm_{m,n}).
$$

Note that by the matrix version of Nehari's theorem, the right-hand side of \rf{dis} is the norm of the Hankel operator $H_\Phi:H^2(\C^n)\to H^2_-(\C^m)$ that is defined on the Hardy class $H^2$ of
$\C^n$-valued functions by
\bay
\label{Han}
H_\Phi f=\pp_-(\Phi f),
\ey
where $\pp_-$ is the orthogonal projection from vector the space $L^2(\C^m)$ onto the subspace
$H^2_-(\C^m)\df L^2(\C^m)\ominus H^2(\C^m)$
(see, e.g.,  \cite{P}, Ch. 2, \S\,2).

Note that this problem is very important in applications in control theory, see e,.g., \cite{F} and 
\cite{P}, Ch. 11.

By the Hahn--Banach theorem,
\bay
\label{HB}
\dist_{L^\be}\big(\Phi, H^\be(\mm_{m,n})\big)=\sup
\left|\int_\T\trace\big(\Phi(\z)\Psi(\z)\big)\,d\m(\z)\right|,
\ey
where the supremum is taken over all matrix functions $\Psi\in H^1_0(\mm_{n,m})$ such that
$\|\Psi\|_{L^1(\bS_1)}=1$. Here $H^1_0(\mm_{n,m})$ is the subspace of the Hardy class
$H^1(\mm_{n,m})$ of $m\times n$ matrix functions vanishing at the origin and the norm $\|A\|_{\bS_1}$
of a matrix $A$ is its {\it trace norm}: $\|A\|_{\bS_1}\df\trace (A^*A)^{1/2}$.

However, it is well known that the infimum is not necessarily attained even for scalar matrix functions (see the Remark following Theorem \ref{123}).
If there exists a matrix function $\Psi\in H^1_0(\mm_{n,m})$ such that
\bay
\label{Up}
\|\Psi\|_{L^1(\bS_1)}=1\quad \mbox{and}\quad
\dist_{L^\be}\big(\Phi, H^\be(\mm_{m,n})\big)=
\int_\T\trace\big(\Phi(\z)\Psi(\z)\big)\,d\m(\z),
\ey
$\Psi$ is called a {\it dual extremal function} of $\Phi$.

Note that the technique of dual extremal functions was used in \cite{Kh} to study the problem of best analytic approximation in the scalar case.

In this paper we characterize the class of matrix functions $\Phi$ that have dual extremal functions.
It turns out that this is equivalent to the fact that the Hankel operator $H_\Phi$ defined by \rf{Han}
has a maximizing vector in $H^2(\C^n)$ which in turn is equivalent to the fact that the matrix function 
$\Phi-F$ (where $F$ is a best approximant to $\Phi$) admits a certain special factorization in terms
thematic matrix functions. The main result will be established in \S\,3. 

In \S\,2 we state Sarason's factorization theorem \cite{S} which will be used in \S\,3 and we define the notion of a thematic matrix function.

\

\section{\bf Preliminaries}
\setcounter{equation}{0}

\

{\bf 1. Sarason's Theorem.} 
We are going to use the following result by D. Sarason:

\medskip

{\bf Sarason's Theorem [S].} {\it Let $\h$ be a separable Hilbert space and 
let $\Psi$ be an analytic integrable $\B(\h)$-valued function on $\T$.
Then there exist analytic square integrable functions $Q$ and $R$
such that 
\bay
\label{sar}
\Psi=QR,\quad R^*R=\big(\Psi^*\Psi\big)^{1/2},\quad
\mbox{and}\quad Q^*Q=RR^*\quad\mbox{a.e. on}~~\T.
\ey}

\medskip

Sarason's theorem implies the following fact: 

\medskip

{\it Let $\Psi$ be a matrix function in $H^1_0(\mm_{n,n})$. Then there
exist matrix functions 
$Q\in H^2(\mm_{n,n})$ and
$R\in H^2_0(\mm_{n,n})$
such that}
$$
\Psi=QR\quad\mbox{and}\quad\|\Psi\|_{L^{1}(\bS_1)}=
\|Q\|_{L^2(\bS_2)}\|R\|_{L^2(\bS_2)}.
$$
Here $H^2_0(\mm_{n,n})$ is the Hardy class of $n\times n$ matrix functions vanishing at the origin.
Recall that the {\it Hilbert--Schmidt norm}  $\|A\|_{\bS_2}$ of a matrix $A$ is defined by 
$\|A\|_{\bS_2}=\trace A^*A$.

\medskip

{\bf 2.  Thematic matrix functions.} The notion of a thematic matrix function was introduced in \cite{PY}.
It turned out that it is very useful in the study of best approximation by analytic matrix functions
(see \cite{P}, Ch. 14). 

Recall that a bounded analytic matrix function $\Theta$ is called an {\it inner function} if 
$\Theta^*(\z)^*\Theta(\z)=I$ form almost all $\z\in\T$, where $I$ is the identical matrix.
A matrix function $F\in H^\be(m,n)$ is called {\it outer} if the operator of multiplication by $F$ on 
$H^2(\C^n)$ has dense range in $H^2(\C^m)$. Finally, we say that a bounded analytic matrix function
$G$ is called {\it co-outer} if the transposed matrix function $G^{\rm t}$ is outer.

An $n\times n$ matrix function $V$ is called a {\it thematic matrix function} if it has the form
$$
V=
\left(
\begin{array}{cc}
\bs{v}  &  \ov{\Theta} 
\end{array}
\right),
$$
where $\bs{v}$ is a column function, both functions $\bs{v}$ and $\Theta$ are inner and co-outer
bounded analytic functions such that $V$ takes unitary values on $\T$, i.e.,
$$
V^*(\z)V(\z)=I,\quad\mbox{for almost all}\quad\z\in\T.
$$
Note that a bounded analytic column function is co-outer if and only if its entries are coprime,
i.e., they do not have a common nonconstant inner factor.

\

\section{\bf The main result}
\setcounter{equation}{0}

\

It is easy to see that a matrix function $\Phi\in L^\be(\mm_{m,n})$ has a dual extremal function
if and only if $\Phi-F$ has a dual extremal function for any $F\in H^\be(\mm_{m,n})$. Moreover,
if $\Psi$ is a dual extremal function of $\Phi$, than $\Psi$ is also a dual extremal function for $\Psi-F$
 for any $F\in H^\be(\mm_{m,n})$. Thus to characterize the class of matrix functions that possess 
 extremal functions, it suffices to consider badly approximable matrix functions.

\begin{thm}
\label{123}
Let $\Phi$ be a nonzero badly approximable function in $L^\be(\mm_{m,n})$
with $m\ge2$ and $n\ge2$.
The following are equivalent:

{\em(i)} the Hankel operator $H_\Phi$ has a maximizing vector;

{\em(ii)} $\Phi$ has a dual extremal function $\Psi\in H^1_0(\mm_{n,m})$;

{\em(iii)} $\Phi$ has a dual extremal function $\Psi\in H^1_0(\mm_{n,m})$ such that
$\rank\Psi(\z)=1$ almost everywhere on $\T$;

{\em(iv)} $\Phi$ admits a factorization
\bay
\label{pt}
\Phi=
W^*\left(\begin{array}{cc}tu&\0\\[.2cm]\0&\Phi_\#\end{array}\right)V^*,
\ey
where $t=\|\Phi\|_{L^\be(\mm_{m,n})}$, $V$ and $W^{\rm t}$
are thematic matrix functions, $u$ is a scalar function of the form
$u=\bar z\bar\vt\bar h/h$ for an inner function $\vt$ and an outer
function $h$ in $H^2$, and $\Phi_\#$ is an $(n-1)\times(m-1)$
matrix function such that $\|\Phi_\#(\z)\|\le t$ for almost all $\z\in\T$.
\end{thm}

Note that the proof of the implication (i)$\imp$(iv) is contained in \cite{PY}, see also 
\cite{P}, Ch. 14, Th. 2.2. However, we give here the proof of this implication for completeness.

\medskip

\Pf (ii)$\imp$(i). By adding zero columns or zero rows if necessary, we
may reduce the general case to the case $m=n$. Let $\Psi$ be a matrix
function in $H^1_0(\mm_{n,n})$ that satisfies \rf{Up}.
 By Sarason's theorem, there
exist functions
$Q\in H^2(\mm_{n,n})$ and
$R\in H^2_0(\mm_{n,n})$
such that
$$
\Psi=QR\quad\mbox{and}\quad1=\|\Psi\|_{L^{1}(\bS_1)}=
\|Q\|_{L^2(\bS_2)}\|R\|_{L^2(\bS_2)}.
$$
Let $e_1,\cdots,e_n$ be the standard orthonormal basis in $\C^n$. We have
\begin{align*}
\int_\T\trace\big(\Phi(\z)\Psi(\z)\big)\,d\m(\z)&=
\int_\T\trace\big(\Phi(\z)Q(\z)R(\z)\big)\,d\m(\z)\\[.2cm]
&=\int_\T\trace\big(R(\z)\Phi(\z)Q(\z)\big)\,d\m(\z)\\[.2cm]
&=\sum_{j=1}^k\int_\T
\big(\Phi(\z)Q(\z)e_j,R^*(\z)e_j\big)\,d\m(\z)\\[.2cm]
&=\sum_{j=1}^k(H_\Phi Qe_j,R^*e_j)
\end{align*}
(we consider here $Q e_j$ and $R^*e_j$ as vector functions).
By the Cauchy--Bunyakovsky--Schwarz inequality, we have
\begin{align*}
\left|\int_\T\trace\big(\Phi(\z)\Psi(\z)\big)\,d\m(\z)\right|&\le
\sum_{j=1}^n|(H_\Phi Qe_j,R^*e_j)|\\[.2cm]
&\le\left(\sum_{j=1}^n\|H_\Phi Qe_j\|_{L^2(\C^n)}^2\right)^{1/2}
\left(\sum_{j=1}^n\|R^*e_j\|_{L^2(\C^n)}^2\right)^{1/2}\\[.2cm]
&
\le\|H_\Phi\|\!
\left(\sum_{j=1}^n\|Qe_j\|_{L^2(\C^n)}^2\!\right)^{1/2}\!\!
\left(\sum_{j=1}^n\|R^*e_j\|_{L^2(\C^n)}^2\!\right)^{1/2}.
\end{align*}
Clearly,
\begin{align*}
\left(\sum_{j=1}^n\|Qe_j\|_{L^2(\C^n)}^2\right)^{1/2}&=
\left(\sum_{j=1}^n\int_\T\|Q(\z)e_j\|_{\C^n}^2\,d\m(\z)\right)^{1/2}\\[.2cm]
&=\left(\int_\T\|Q(\z)\|_{\bS_2}^2\,d\m(\z)\right)^{1/2}
=\|Q\|_{L^2(\bS_2)}.
\end{align*}
and
$$
\left(\sum_{j=1}^n\|R^*e_j\|_{L^2(\C^n)}^2\!\right)^{1/2}=
\|R\|_{L^2(\bS_2)}.
$$
Since $\Phi$ is badly approximable, 
we have $\|H_\Phi\|=\|\Phi\|_{L^\be(\mm_{n,n})}$.

It follows that
\begin{align*}
\|\Phi\|_{L^\be(\mm_{n,n})}&=
\left|\int_\T\trace\big(\Phi(\z)\U(\z)\big)\,d\m(\z)\right|\\[.2cm]
&\le
\left(\sum_{j=1}^n\|H_\Phi Qe_j\|_{L^2(\C^n)}^2\right)^{1/2}
\left(\sum_{j=1}^n\|R^*e_j\|_{L^2(\C^n)}^2\right)^{1/2}\\[.2cm]
&
\le\|H_\Phi\|\cdot
\|Q\|_{L^2(\bS_2)}\|R\|_{L^2(\bS_2)}=\|\Phi\|_{L^\be(\mm_{n,n})}.
\end{align*}
Thus all inequalities are equalities and if $Qe_j\ne\0$, then
$Qe_j$ is a maximizing vector of $H_\Phi$.

The implication (iii)$\imp$(ii) is trivial.

(iv)$\imp$(iii). Suppose that $\Phi$ is a function given by
\rf{pt}. By multiplying $h$ by a constant if necessary,
we may assume without loss of generality that $\|h\|_{L^2}=1$. Let
$$
V=\left(\begin{matrix}\bs{v}&\ov{\Theta}\end{matrix}\right)\quad{and}
\quad
W^{\rm t}=\left(\begin{matrix}\bs{w}&\ov{\Xi}\end{matrix}\right).
$$

Put
$$
\Psi=z\vt h^2\left(\begin{matrix}\bs{v}&\0\end{matrix}\right)
\left(\begin{matrix}\bs{w}^{\rm t}\\[.2cm]\0\end{matrix}\right).
$$
Clearly,
$$
\|\Psi\|_{L^1(\bS_1)}=\|h^2\|_{L^1}=1
$$
and it is easy to see that
\begin{align*}
\int_\T\trace\big(\Phi(\z)\Psi(\z)\big)\,d\m(\z)&=
\int_\T z\vt h^2\trace\left(\left(\begin{matrix}\bs{w}^{\rm t}
\\[.2cm]\0\end{matrix}\right)\Phi
\left(\begin{matrix}\bs{v}&\0\end{matrix}\right)\right)\,d\m\\[.2cm]
&=\int_\T\trace\left(\begin{matrix}|h|^2&\0\\[.2cm]\0&\0\end{matrix}\right)
\,d\m=1.
\end{align*}

(i)$\imp$(iv). Let $f$ be a maximizing vector of $H_\Phi$. It is well known (see \cite{P}, Ch. 2. Th. 2.3)
that 
$$
\|\Phi(\z)\|_{\mm_{m,n}}=\|\Phi\|_{L^\be}=\|H_\Phi\|,\quad
\|\Phi(\z)f(\z)\|_{\C^m}=\|H_ \Phi\|\cdot\|f(\z)\|_{\C^n},\quad\quad\z\in\T,
$$
and
$$
\Phi f\in H^2_-(\C^m),.
$$
Put
$$
g=\frac1{\|H_\Phi\|}\bar z\ov{\Phi f}=\frac1{\|H_\Phi\|}\bar z\ov{H_\Phi f}\in H^2(\C^m).
$$
Then 
$$
\|f(\z)\|_{\C^n}=\|g(\z)\|_{\C^m},\quad\z\in\T.
$$
It follows that both $f$ and $g$ admit factorizations
$$
f=\vt_1h\bs{v},\quad g=\vt_2h\bs{w},
$$
where $\vt_1$ and $\vt_2$ are scalar inner functions, $h$ is a scalar outer function in $H^2$, and
$\bs{v}$ and $\bs{w}$ are inner and co-outer column functions.
Then $\bs{v}$ and $\bs{w}$ admit thematic completions, i.e., there are inner and co-outer
matrix functions $\Theta$ and $\Xi$ such that the matrix functions
$$
\left(
\begin{array}{cc}
\bs{v}  &  \ov{\Theta} 
\end{array}
\right)
\quad\mbox{and}\quad
\left(
\begin{array}{cc}
\bs{w}  &  \ov{\Xi} 
\end{array}
\right)
$$
are thematic. Put
$$
V=\left(
\begin{array}{cc}
\bs{v}  &  \ov{\Theta} 
\end{array}
\right),
\quad
W=\left(
\begin{array}{cc}
\bs{w}  &  \ov{\Xi} 
\end{array}
\right)^{\rm t},\quad\mbox{and}\quad u=\bar z\ov{\vt}_1\ov{\vt}_2\bar h/h.
$$
Consider the matrix function $W\Phi V$. It is easy to see that its upper left entry is equal to
$$
\bs{w}^{\rm t}\Phi\bs{v}=\frac{\ov{\vt}_2}{h}g^{\rm t}\Phi\frac{\ov{\vt}_1}{h}f
=\|H_\Phi\|\bar z\frac{\ov{\vt}_1\ov{\vt}_2}{h^2}g^{\rm t}g=
\|H_\Phi\|u=tu.
$$

Since  the norm of $(W\Phi V)(\z)$ is equal to $t$  and its upper left entry
$tu(\z)$ has modulus $t$ almost everywhere, it is easy to see that the matrix function
$W\Phi V$ has the form
$$
W\Phi V=
\left(
\begin{array}{cc}
 tu &\0      \\
\0  & \Phi_\#    
\end{array}
\right),
$$
where $\Phi_\#$ is an $(m-1)\times(n-1)$ matrix function such that $\|\Phi_\#\|_{L^\be}\le t$. It follows that
$$
\Phi=W^*
\left(
\begin{array}{cc}
 tu &\0      \\
\0  & \Phi_\#    
\end{array}
\right)V^*
$$
which completes the proof. $\bl$

\medskip

{\bf Remark.}  In the case when $\Phi$ has size $m\times1$, $m>1$, Theorem \ref{123} remains true
if we replace the factorization in \rf{pt} with the factorization
$$
\Phi=W^*\left(\begin{array}{c}tu\\\0\end{array}\right)
$$
where $W^{\rm t}$ is a thematic matrix function and $u$ has the form $u=\bar z\bar\vt\bar h/h$, where
$\vt$ is a scalar inner function and $h$ is an scalar outer function in $H^2$.

Similarly, the theorem can be stated in the case of size $1\times n$, $n>1$.

In the case of scalar functions, the result also holds if we replace (iv) with the condition that $\Phi$
admits a factorization in the form
$$
\Phi=\bar z\bar\vt\bar h/h,
$$
where $\vt$ is a scalar inner function and $h$ is an scalar outer function in $H^2$.

\medskip

Since it is well known that not all scalar badly approximable functions have constant modulus on $\T$
(see e.g., \cite{P}, Ch. 1, \S\,1), there are scalar functions in $L^\be$ that have no dual
extremal functions.

\

\

\noindent
\begin{tabular}{p{8cm}p{14cm}}
V.V. Peller \\
Department of Mathematics \\
Michigan State University  \\
East Lansing, Michigan 48824\\
USA
\end{tabular}


\begin{thebibliography}{99}

\bibitem[F]{F} {\sc B. A. Francis}, {\it A Course in $H^\infty$ Control Theory,} 
Lecture Notes in
Control and Information Sciences No. 88, Springer Verlag, Berlin, 1986.
\bibitem[Kh]{Kh} {\sc S. Khavinson}, On some extremal problems of the theory 
of analytic functions, {\em Uchen. Zapiski Mosk. Universiteta, Matem.} {\bf 144:4}
(1951), 133--143. English Translation: Amer. Math. Soc. Translations (2) {\bf 32}
(1963), 139--154.
\bibitem[P]{P} {\sc V.V. Peller}, {\em Hankel operators and their applications,}
 Springer-Verlag, New York, 2003.
\bibitem[PY]{PY} {\sc V.V. Peller and N.J. Young}, Superoptimal analytic 
approximations of 
matrix functions, {\em J. Funct. Anal.} {\bf 120} (1994), 300-343. 
\bibitem[S]{S} {\sc D. Sarason}, {\it Generalized interpolation in $H^\infty$,} Trans.
Amer. Math. Soc., {\bf 127} (1967) \lb179--203.


\end{thebibliography}
\end{document}